\documentclass[11pt,a4paper,fleqn]{article}
\usepackage{a4wide,amsfonts,amsmath,latexsym,amssymb,euscript,graphicx,units,mathrsfs}

\usepackage{graphicx}
\usepackage{color}
\usepackage{amssymb}
\usepackage{amssymb}
\usepackage[T1]{fontenc}
\usepackage{latexsym}
\usepackage{xypic}
\usepackage{eufrak}
\usepackage{euscript}
\usepackage{amsfonts,amsmath}
\usepackage{verbatim}
\usepackage{fancyhdr}
\usepackage[english]{babel}
\usepackage{mathrsfs}
\usepackage{units}

\newtheorem{prop}{Proposition}[section]
\newtheorem{cor}[prop]{Corollary}
\newtheorem{lemma}[prop]{Lemma}

\newtheorem{theorem}[prop]{Theorem}
\newtheorem{defi}[prop]{Definition}

\renewcommand{\geq}{\geqslant}
\def\leq{\leqslant}

\newcommand{\N}{\mathbb{N}}
\newcommand{\Z}{\mathbb{Z}}
\newcommand{\R}{\mathbb{R}}

\def \P{\mathbb{P}}

\def\HH{\EuFrak H}

\def\e{\varepsilon}

\def\1{{\mathbf{1}}}

\def\1{{\mathbf{1}}}
\def\0.5{{\frac{1}{2}}}

\def\E{\mathbb{E}}

\newcommand{\qed}{\nopagebreak\hspace*{\fill}
{\vrule width6pt height6ptdepth0pt}\par}

\date{ }
\begin{document}
\title{ {\bf Almost sure limit theorems on Wiener chaos:\\
the non-central case } }

\author{Ehsan Azmoodeh \\ \small{Ruhr-Universit\"at Bochum, Fakult\"at f\"ur Mathematik} \\ \small{{ \tt ehsan.azmoodeh@rub.de}} \medskip \\
	Ivan Nourdin \\ \small{Universit\'e du Luxembourg, Unit\'e de Recherche en Math\'ematiques} \\ \small{ {\tt ivan.nourdin@uni.lu}}}

\maketitle

\abstract
In \cite{BNT}, a framework to prove
almost sure central limit theorems for sequences $(G_n)$ belonging to the Wiener space was developed, with a particular emphasis of the case where $G_n$ takes the form of a multiple Wiener-It\^o integral with respect to a given isonormal Gaussian process.
In the present paper, we complement the study initiated in \cite{BNT},  by considering the more general situation where the sequence $(G_n)$ may not need to converge to a Gaussian distribution.
As an application, we prove that partial sums of Hermite polynomials of increments
of fractional Brownian motion satisfy an almost sure limit theorem in the long-range dependence case, a problem left open in \cite{BNT}.

\vskip0.3cm
\noindent {\bf Keywords}: Almost sure limit theorem; Multiple Wiener-It\^o integrals; Malliavin Calculus; Characteristic function; Wiener chaos; Hermite distribution; Fractional Brownian motion.\\

\noindent{\bf MSC 2010}: {60F05, 60G22, 60G15, 60H05, 60H07}.

\section{Introduction and main results}

Let us start with the following definition, which plays a pivotal role in the paper.
\begin{defi}\label{def:ASLT}
Let $(G_n)$ be a sequence of real-valued random variables defined on a common probability space $(\Omega,\mathcal{F},\mathbb{P})$ and let $\mu$ be a probability measure on $(\R,\mathcal{B}(\R))$.
\begin{enumerate}
\item We say that  $(G_n)$ satisfies an \underline{almost sure limit theorem} (in short: \underline{ASLT}) with
limit $\mu$ if
\begin{equation}\label{defi1}
\P \left( \frac{1}{\log n} \sum_{k \le n} \frac{1}{k}  \delta_{G_k} \Rightarrow\mu \mbox{ as $n\to\infty$}\right)=1,
\end{equation}
where $\delta_x$ denotes the Dirac point measure of $x$ and `$\Rightarrow$' stands for the weak convergence of measures.
\item If (\ref{defi1}) holds true with $\mu= N(0,1)$ (the standard Gaussian distribution),
we say that $(G_n)$ satisfies an \underline{almost sure central limit theorem} (in short: \underline{ASCLT}).
\end{enumerate}
\end{defi}
When $\mu$ has a density, it follows from a classical approximation argument  that (\ref{defi1}) is equivalent to saying that
\begin{equation}\label{defi1bis}
\P\left(\lim_{n\to\infty}\frac{1}{\log n} \sum_{k \le n} \frac{1}{k}  {\bf 1}_{\{ G_k\leq x\}} = \mu(-\infty,x]\mbox{ for all $x\in\R$}\right)=1.
\end{equation}

Convergence (\ref{defi1}) is a kind of weighted \emph{strong} law of large numbers, but related to a quantity $G_n$ converging \emph{in law} to $\mu$; hence its name. 
Historically, the first appearance of an ASLT
was in the book \cite{Levy} by L\'evy. It was stated without proof in the simple situation where $\mu= N(0,1)$ and $G_n$ has the form $G_n=(X_1+\ldots+X_n)/\sqrt{n}$, for i.i.d. \!summands $X_i$ with mean 0 and variance 1.
L\'evy's ASCLT went unnoticed for more than a half century,  until it was exploited by 
Brosamler \cite{Bro} and Schatte \cite{Sch} independently in 1998.
Since then, there has been a lot of activities around the ASLTs in various contexts, see e.g. \cite{BC,BD,Gon,IL98,IL00,PS} and references therein.

We stress on the fact that the validity (or not) of (\ref{defi1}) is a property of
the distribution of the \emph{whole} sequence $(G_n)$. Let us elaborate on this point. If $G_n\Rightarrow\mu$, it is well-known by Skorokhod's representation theorem 
that there exists a family $(G_n^\star)_{n\in\N\cup\{\infty\}}$ of random variables defined on a common probability space such that
$G_n^\star\overset{\rm law}{=}G_n$ for all $n$ and $G_n^\star\to G_\infty^\star\sim \mu$ almost surely. Moreover, Ces\`aro summation theorem implies that, almost surely,
$(\log n)^{-1} \sum_{k \le n} \varphi( G^\star_k)/k  \longrightarrow \varphi(G^\star_\infty)$ for any continuous bounded function $\varphi:\R\to\R$; that is, $(G_n^\star)$ does \emph{not} satisfy an ASLT except if $\mu$ is a Dirac mass, since apart in this case $\varphi(G^\star_\infty)\neq 
\E[\varphi(G^\star_\infty)]$ for at least one function $\varphi$.

Now we have made 
these preliminary comments on general ASLTs, let us 
concentrate on the specific kind of random variables $(G_n)$ 
we are interested in in this paper: namely, random variables taking the form of \emph{multiple Wiener-It\^o integrals}. 
Since the discovery by Nualart and Peccati \cite{NuPe} of their fourth moment theorem (according to which 
a normalized sequence of multiple Wiener-It\^o integrals  converges in law to $N(0,1)$ if and only if its fourth moment converges to 3),
those stochastic integrals have become a probabilistic object of considerable interest. See for instance the 
constantly updated webpage \cite{webpage} for
a demonstration of the intense activity surrounding them.

More specifically, a  sequence $(G_n)$  of multiple Wiener-It\^o integrals converging in law to some $\mu$ being given, our goal in this paper is to provide a meaningful set of conditions under which an ASLT holds true.
The central case (that is, when $\mu = N(0,1)$)
has already been the object of \cite{BNT}. Therein, 
a framework has been developed in order to show an ASCLT, 
by taking advantage of the numerous estimates coming from the Malliavin-Stein approach \cite{n-p-book}.
One of the main results of \cite{BNT}, that we state here for the sake of comparison, is the following. The notation $I_q(g)$ refers to the $q$th Wiener-It\^o integral with kernel $g$ associated to a given isonormal Gaussian process $X$ over a separable Hilbert space  $\HH$, whereas the notation $\otimes_r$ stands for the contraction operator, see Section 2.3 for more details.

\begin{theorem}\label{thm1}
Fix $q\geq 2$, and let $(G_n)$ be a sequence of the form $G_n=I_q(g_n)$, 
$n\geq 1$. Suppose (for simplicity) that $\E[G_n^2]=1$ for all $n$.
Suppose in addition that the following two conditions are met.
\begin{enumerate}
\item[(A1)] for every $r\in\{1,\ldots,q-1\}$: $\sum_{n\geq 2}\frac{1}{n\log^2n}\sum_{k=1}^n\frac1k\|g_k\otimes_r g_k\|_{\HH^{\otimes (2q-2r)}}<\infty$;
\item[(A2)] $\sum_{n\geq 2}\frac{1}{n\log^3n}\sum_{k,l=1}^n \frac{|\E[G_kG_l]|}{kl}<\infty$.
\end{enumerate}
Then $(G_n)$ satisfies an ASCLT, that is, (\ref{defi1}) holds true with 
$\mu = N(0,1)$.
\end{theorem}

An application of high interest of Theorem \ref{thm1} studied in \cite{BNT} concerns
the celebrated Breuer-Major theorem.
Let $\{X_k\}$ be a stationary centered Gaussian family, characterized by 
its correlation function $\rho:\Z\to\R$ given by $\rho(k-l):=\E[X_kX_l]$, and assume that $\rho(0)=1$.
Fix $q\geq 2$ and let $H_q$ denote the $q$th Hermite polynomial. 
Finally, consider the sequence of the Hermite variations of $X$,  defined as
\begin{equation}\label{vn}
V_n=\sum_{k=1}^n H_q(X_k),
\quad n\geq 1.
\end{equation}
It is well known since the seminal works by Breuer and Major \cite{BM}, Giraitis and Surgailis \cite{GS} and Taqqu \cite{T}, that the study of the fluctuations of $V_n$ crucially depends on the summability of $|\rho|^q$. More precisely, if we
set
\begin{equation}\label{gn}
G_n=\frac{V_n}{\sqrt{{\bf Var}(V_n)}},\quad n\geq 1,
\end{equation}
then $G_n\Rightarrow N(0,1)$ as soon as
$\sum_{k\in\Z} |\rho(k)|^q<\infty$.
When $\sum_{k\in\Z} |\rho(k)|^q=\infty$, the sequence $(G_n)$ may
converge to a Gaussian (in some critical situations) or \emph{non}-Gaussian limit (more likely), see below.

A classical example of a stationary Gaussian sequence $\{X_k\}_{k\in\Z}$ falling within this framework is given by the fractional Gaussian noise
\begin{equation}\label{eq:fgn}
X_k=B_{k+1}-B_{k}, 
\end{equation}
where $B=(B_t)_{t \in \R}$ is a fractional Brownian motion of Hurst index $H\in (0,1)$, that is, $B$ is  a centered Gaussian process with covariance function $$ \E[B_t B_s] := \frac{1}{2} \left\{  \vert t \vert ^{2H} + \vert s \vert ^{2H} - \vert t-s \vert^{2H}   \right\}, \quad s,t \in \R.$$
More precisely, it is nowadays well-known that
\begin{itemize}
\item[(i)] if $0<H\leq 1-\frac1{2q}$ then $G_n\Rightarrow N(0,1)$;
\item[(ii)] if $1-\frac1{2q}<H<1$ then $G_n\Rightarrow\mu_{H,q}$,  the Hermite distribution with parameters $H$ and $q$.
\end{itemize}
The proof of (i) follows directly from the Breuer-Major theorem when $H<1-\frac1{2q}$, whereas a little more effort are required in the critical case $H=1-\frac{1}{2q}$ (see, e.g., \cite{BN}). For a definition of the Hermite distribution $\mu_{H,q}$ as well as a short proof of the weak convergence (ii) the reader can consult, e.g., \cite[Proposition 6.1]{BNT} and \cite[Proposition 7.4.2]{n-p-book}. \\

In this paper, we aim to answer the following question: \emph{can we associate an almost sure limit theorem to the previous two convergences in law (i) and (ii)?} This problem is actually not new, and has been first considered almost one decade ago.
In \cite{BNT}, a positive answer has been indeed given for (i),
with the help of Theorem \ref{thm1}:
namely, if $H\leq 1-\frac1{2q}$ and if $G_n$ is defined by (\ref{gn}), then 
(\ref{defi1}) holds true with $\mu= N(0,1)$.
What about the case $H\in(1-\frac1{2q},1)$?
A solution to this problem was left open in \cite{BNT} probably because, so far,
the Malliavin-Stein approach used therein 
has been mainly developed
in the case of {\it normal} approximation.

In the present paper, our main motivation is then to study the case $H\in(1-\frac1{2q},1)$
left open in \cite{BNT}.
In order to do so, we first develop a general abstract result (Theorem \ref{thm2}), of independent interest and valid in a framework similar to that of Theorem \ref{thm1};
then, we apply it to our specific situation (Corollary \ref{thm3}).

\begin{theorem}\label{thm2}
Fix $q\geq 2$, and let $(G_n)$ be any sequence of the form $G_n=I_q(g_n)$, 
$n\geq 1$. Suppose $G_n\Rightarrow\mu$, and
denote by $\phi_n$ (resp. $\phi_\infty$) the characteristic function of $G_n$ (resp.
$\mu)$.
Suppose in addition that the following two conditions are met.
\begin{enumerate}
\item[(B1)] for all $r>0$: $\sup_{|t|\leq r}\sum_{n\geq 2}\frac{1}{n\log^3n}\left|\sum_{k=1}^n\frac1k(\phi_k(t)-\phi_\infty(t))\right|^2<\infty$;
\item[(B2)] $\sum_{n\geq 2}\frac{1}{n\log^3n}{\bf Var}\left(\sum_{k=1}^n \frac{G_k^2}{k}\right)<\infty$.
\end{enumerate}
Then $(G_n)$ satisfies an ASLT, that is, (\ref{defi1}) holds true.
\end{theorem}

Let $G_n$ be given by (\ref{gn}) with $H$ belonging to
$(1-\frac{1}{2q},1)$.
In \cite[Section 6]{BNT}, the authors have constructed an \emph{explicit} sequence $(\widetilde{G}_n)$ such that 
$\widetilde{G}_n\overset{\rm law}{=}G_n$ for any $n$ and
$\widetilde{G}_n$ does \emph{not} satisfy an ASLT.
But as we already explained in the second paragraph following Definition \ref{defi1}, from this there is nothing to learn about the validity or not of an ASLT associated to the original sequence $(G_n)$. 
Here, using our abstract Theorem \ref{thm2}, we are actually able to prove that $(G_n)$ \emph{does} satisfy an ASLT.

\begin{cor}\label{thm3}
Let $G_n$ be given by (\ref{gn}), with $X_k$ as in \eqref{eq:fgn} and $H\in(1-\frac{1}{2q},1)$.
Then $(G_n)$ satisfies an ASLT with $\mu=\mu_{H,q}$  the Hermite distribution with parameters $H$ and $q$.
\end{cor}

To conclude this introduction, we would like to stress that
our Theorem \ref{thm2} is a true extension of Theorem \ref{thm1}, in the sense that the latter can be obtained as a particular case of the former: see Section 5 for the details.

Also, for the sake of illustration, we develop 
in Section \ref{sec-easyexample} an 
easy example of application of our Theorem \ref{thm2} involving 
a sequence of independent standard Gaussian random variables.

The rest of the paper is organized as follows.
In Section \ref{preliminaries} several preliminary results that are needed for the proofs are collected.
Section \ref{section3} is devoted to the proof of Theorem \ref{thm2}.
In Section \ref{sec-easyexample}, an easy application of Theorem \ref{thm2}
to partial sums of independent standard Gaussian random variables is presented.
The proof of Corollary \ref{thm3} can be found in Section \ref{section 4}.
Finally, Section \ref{section 5} details how to deduce 
\cite[Theorem \ref{thm1}]{BNT}
from our Theorem \ref{thm2}.

\section{Preliminaries}\label{preliminaries}

\subsection{Ibragimov-Lifshits criterion for ASLT}\label{sec:IL}

As we will see, we are not going to prove \emph{directly} that  $(G_n)$ in Theorem \ref{thm2} satisfies (\ref{thm1}). Instead, we are going to check
the validity of 
a powerful criterion due to Ibragimov and Lifshits \cite{IL00}, that we recall here for the convenience of the reader.

\begin{theorem}\label{thm:generalAS}
	Let $(G_n)$ be a sequence of real-valued random variables defined on a common probability space and converging in distribution towards  $\mu$, and set 
	\begin{equation}\label{eq:Delta}
	\Delta_n(t) = \frac{1}{\log n} \sum_{k=1}^{n} \frac{1}{k} \left( e^{itG_k} - \int_\R e^{itx}\mu(dx) \right).
	\end{equation}
	If
	\begin{equation}\label{IL-condition}
	\sup_{\vert t \vert \le r} \sum_{n} \frac{\E \vert \Delta_n(t) \vert^2}{n \log n} < \infty\quad\mbox{for all $r >0$},
	\end{equation}
	then $(G_n)$ satisfies the ASLT (\ref{defi1}).
\end{theorem} 

\subsection{An easy reduction lemma}

In this section, we state and prove an easy reduction lemma, that we are going to
use in the proof of Corollary \ref{thm3}.

\begin{lemma}\label{reduction}
Suppose that $(G_n)$ is such that $G_n\Rightarrow\mu$,
and assume that $\mu$ has a density. Let $(a_n)$ be a real-valued sequence converging to $a_\infty\neq 0$.
Then $(G_n)$ satisfies an ASLT if and only if $(a_nG_n)$ does.
\end{lemma}
{\it Proof}. 
Without loss of generality, we may and will assume that $a_n>0$ for all $n$ and that $a_\infty=1$.
By symmetry, only the implication ``if $(G_n)$ satisfies an ASLT then $(a_nG_n)$ does'' has to be proved. 

So, let us assume that $(G_n)$ satisfies an ASLT. Since $\mu$ has a density, we are going to use criterion (\ref{defi1bis}).
Fix $\e>0$, and let $k_0$ be such that $1-\e\leq a_k\leq 1+\e$ for all $k\geq k_0$.
For any $x\in\R$, we can write, for all $k\geq k_0$,
\begin{eqnarray*}
{\bf 1}_{\{G_k\leq\frac{x}{1+\e}\}}{\bf 1}_{\{x\geq 0\}}
+{\bf 1}_{\{G_k\leq\frac{x}{1-\e}\}}{\bf 1}_{\{x<0\}}
&\leq& 
{\bf 1}_{\{a_kG_k\leq x\}}\\
&\leq& 
{\bf 1}_{\{G_k\leq\frac{x}{1-\e}\}}{\bf 1}_{\{x\geq 0\}}
+{\bf 1}_{\{G_k\leq\frac{x}{1+\e}\}}{\bf 1}_{\{x<0\}}.
\end{eqnarray*}
As a consequence, 
$$
\frac{1}{\log n}\sum_{k=k_0}^n \frac{1}{k}\left(
{\bf 1}_{\{G_k\leq\frac{x}{1+\e}\}}{\bf 1}_{\{x\geq 0\}}
+{\bf 1}_{\{G_k\leq\frac{x}{1-\e}\}}{\bf 1}_{\{x<0\}}
\right)
\leq 
 \frac{1}{\log n}\sum_{k=1}^n \frac{1}{k}{\bf 1}_{\{a_kG_k\leq x\}}
 $$
 and
 $$
 \frac{1}{\log n}\sum_{k=1}^n \frac{1}{k}{\bf 1}_{\{a_kG_k\leq x\}}
\leq
\frac{k_0}{\log n}
+\frac{1}{\log n} \sum_{k=k_0}^n \frac{1}{k}\left({\bf 1}_{\{G_k\leq\frac{x}{1-\e}\}}{\bf 1}_{\{x\geq 0\}}
+{\bf 1}_{\{G_k\leq\frac{x}{1+\e}\}}{\bf 1}_{\{x<0\}}\right).
$$
Letting $n\to\infty$, we deduce from (\ref{defi1bis}) that, almost surely,
for all $x\in\R$,
\begin{multline*}
\mu(-\infty,(1+\e)^{-1}x]\,{\bf 1}_{\{x\geq 0\}}
+\mu(-\infty,(1-\e)^{-1}x]\,{\bf 1}_{\{x< 0\}}
\leq 
\liminf_{n\to\infty}
\frac{1}{\log n} \sum_{k=1}^n \frac{1}{k}{\bf 1}_{\{a_kG_k\leq x\}}\\
\leq 
\limsup_{n\to\infty}
\frac{1}{\log n} \sum_{k=1}^n \frac{1}{k}{\bf 1}_{\{a_kG_k\leq x\}}
\leq 
\mu(-\infty,(1-\e)^{-1}x]\,{\bf 1}_{\{x\geq 0\}}
+\mu(-\infty,(1+\e)^{-1}x]\,{\bf 1}_{\{x< 0\}}.
\end{multline*}
Letting now $\e\to 0$ yields 
$\frac{1}{\log n}\sum_{k=1}^n \frac{1}{k}{\bf 1}_{\{a_kG_k\leq x\}}\to 
\mu(-\infty,x]$,
which is desired conclusion.\qed

\subsection{Elements of Malliavin calculus}\label{sec:EMC}
Our proofs of Theorem \ref{thm2} and Corollary \ref{thm3} are based
on the use of the Malliavin calculus. This is why in this section 
we recall the few elements of Gaussian analysis and Malliavin calculus that will be needed. 
For more details or missing proofs, we invite the reader to 
consult one of the three books \cite{n-p-book, Nbook,DavidEulaliaBook}.

Let $ \HH$ be a real separable Hilbert space. For any $q\geq 1$, we write $ \HH^{\otimes q}$ and $ \HH^{\odot q}$ to
indicate, respectively, the $q$th tensor power and the $q$th symmetric tensor power of $\HH$.

 We denote by $X=\{X(h) : h\in  \HH\}$
an \emph{isonormal Gaussian process} over $ \HH$, that is, $X$ is a centered Gaussian family defined a common probability space $(\Omega,\mathcal{F},\mathbb{P})$ satisfying
$\E\left[ X(g)X(h)\right] =\langle g,h\rangle _{ \HH}$. We also assume that $\mathcal{F}$ is generated by $X$, and use the shorthand notation $L^2(\Omega) = L^2(\Omega, \mathcal{F}, \P)$.

For every $q\geq 1$, the $q$th \textit{Wiener chaos} of $X$ is  defined as the closed linear subspace of $L^2(\Omega)$
generated by the family $\{H_{q}(X(h)) : h\in  \HH,\left\| h\right\| _{ \HH}=1\}$, where $H_{q}$ is the $q$th \emph{Hermite polynomial}:
\begin{equation}\label{hq}
H_q(x) = (-1)^q e^{\frac{x^2}{2}}
\frac{d^q}{dx^q} \big( e^{-\frac{x^2}{2}} \big).
\end{equation}
For
any $q\geq 1$, the mapping $I_{q}(h^{\otimes q})=H_{q}(W(h))$ (for $h\in  \HH,\left\| h\right\| _{ \HH}=1$) can be extended to a
linear isometry between the symmetric tensor product $ \HH^{\odot q}$
(equipped with the modified norm $\sqrt{q!}\left\| \cdot \right\| _{ \HH^{\otimes q}}$)
and the $q$th Wiener chaos.

Let $\{e_{k},\,k\geq 1\}$ be a complete orthonormal system in $\HH$. Given $f\in  \HH^{\odot p}$ and $g\in \HH^{\odot q}$, for every
$r=0,\ldots ,p\wedge q$, the \textit{contraction} of $f$ and $g$ of order $r$
is the element of $ \HH^{\otimes (p+q-2r)}$ defined by
\begin{equation}
f\otimes _{r}g=\sum_{i_{1},\ldots ,i_{r}=1}^{\infty }\langle
f,e_{i_{1}}\otimes \ldots \otimes e_{i_{r}}\rangle _{ \HH^{\otimes
		r}}\otimes \langle g,e_{i_{1}}\otimes \ldots \otimes e_{i_{r}}
\rangle_{ \HH^{\otimes r}}.  \label{v2}
\end{equation}
Notice that the definition of $f\otimes_r g$ does not depend
on the particular choice of $\{e_k,\,k\geq 1\}$, and that
$f\otimes _{r}g$ is not necessarily symmetric; we denote its
symmetrization by $f\widetilde{\otimes }_{r}g\in  \HH^{\odot (p+q-2r)}$.
Moreover, $f\otimes _{0}g=f\otimes g$ equals the tensor product of $f$ and
$g$ while, for $p=q$, $f\otimes _{q}g=\langle f,g\rangle _{ \HH^{\otimes q}}$. 
Contractions appear naturally in the \emph{product formula} for multiple integrals: if $f\in  \HH^{\odot p}$ and $g\in  \HH^{\odot q}$, then
\begin{eqnarray}\label{multiplication}
I_p(f) I_q(g) = \sum_{r=0}^{p \wedge q} r! {p \choose r}{ q \choose r} I_{p+q-2r} (f\widetilde{\otimes}_{r}g).
\end{eqnarray}

\section{Proof of Theorem \ref{thm2}}\label{section3}
As we said in Section \ref{sec:IL}, we are going to check 
the condition \eqref{IL-condition} in our specific framework. 
For $t\in\R$, recall the quantity $\Delta_n(t)$ from \eqref{eq:Delta} and
that $\phi_k$ (resp. $\phi_\infty$) stands for the characteristic function of $G_k$ (resp. $\mu$).
We also denote by $\varphi_{k,l}$ the characteristic function of $G_k-G_l$.
 With these notation in mind, we can now write
\begin{equation*}
\begin{split}
&\E \vert \Delta_n (t) \vert^2 \\
& = \frac{1}{\log^2 n} \sum_{k,l=1}^{n} \frac{1}{kl} \E \Bigg( \Big( e^{it G_k} - \phi_{\infty}(t) \Big) \times \Big( e^{-it G_l} - \phi_{\infty} (-t)\Big) \Bigg) \\
&=\frac{1}{\log^2 n} \sum_{k,l=1}^{n} \frac{1}{kl} \Bigg\{ \Big( \phi_{k} (t)- \phi_{\infty} (t)\Big)  \Big( \phi_{l}(-t) - \phi_{\infty}(-t) \Big) +  \Big( \varphi_{k,l} (t) - \phi_{k}(t) \phi_{l}(-t) \Big) \Bigg\}\\
& =: A_1(n,t) + A_2(n,t). 
\end{split}
\end{equation*}
Now, note that 
\[
A_1(n,t) = \frac{1}{\log^2 n}   \Big \vert \sum_{k=1}^{n} \frac{1}{k} \left( \phi_k (t) - \phi_\infty(t) \right) \Big\vert^2,
\] 
so that 
$
\sup_{|t|\leq r}\sum_{n\geq 2}\frac{A_1(n,t)}{n\log n}<\infty
$ for all $r>0$,
by
condition (B1) in Theorem \ref{thm2}. 

On the other hand, \cite[Proposition 3.1]{NNP}
provides the existence of a universal constant $c$
such that
\[
\Big \vert \varphi_{k,l} (t) - \phi_{k}(t) \phi_{l}(-t) \Big \vert
\leq c|t|\, {\bf Cov}(G^2_k,G^2_l) .
\]
We deduce that
\begin{equation}
\begin{split}
A_2(n,t) &\le  \frac{c \vert t \vert}{\log^2 n} \sum_{k,l=1}^{n} \frac{1}{kl} {\bf Cov}(G^2_k,G^2_l) 
=  \frac{c \vert t \vert}{\log^2 n} {\bf Var} \left( \sum_{k=1}^{n} \frac{G^2_k}{k}\right),
\end{split}
\end{equation}
implying in turn by condition (B2) that
$
\sup_{|t|\leq r}\sum_{n\geq 2}\frac{A_2(n,t)}{n\log n}<\infty
$ for all $r>0$,
and concluding the proof.
\qed

\section{An easy illustrating example}\label{sec-easyexample}

Before proceeding with the proof of Corollary \ref{thm3}, let us illustrate
the use of our Theorem \ref{thm2} on a simple example.

Let $X_0,X_1,X_2\ldots\sim N(0,1)$ be a sequence of independent copies
defined on the same probability space, and consider
$$
G_n=\frac{X_n}{\sqrt{2n}}\sum_{j=0}^{n-1} (X_j^2-1).
$$

A first observation is that $G_n$ can be realized as a Wiener-It\^o integral of order 3 with respect to an isonormal Gaussian process $X$ over the Hilbert space $\HH=L^2(\R_+)$. Indeed, if we identify $X_j$ with $X({\bf 1}_{[j,j+1]})$, we get
\[
G_n=I_3(g_n),\quad \mbox{where $g_n=\frac{1}{\sqrt{2n}}\sum_{j=0}^{n-1}\widetilde{{\bf 1}}_{[n,n+1]\times [j,j+1]^2}$},
\]
where the tilde means that the function has been symmetrized.

Thus, we are in a position to make use of Theorem \ref{thm2}, with $\mu$ the distribution of the product of two independent standard Gaussian random variables.

First, the fact that $G_n\Rightarrow\mu$ is obvious by the classical central limit theorem, since $$G_n\overset{\rm law}{=}
\frac{X_0}{\sqrt{2n}}\sum_{j=1}^{n} (X_j^2-1).$$ 
It is moreover clear that, for any (fixed) 
$r>0$
$$
\sup_{|t|\leq r}\left|\E[e^{itG_k}]-\int_\R e^{itx}d\mu\right|
=\sup_{|t|\leq r}\left|\E\left[e^{-\frac{t^2}{2}\left\{\frac1{\sqrt{2k}}\sum_{j=1}^{n} (X_j^2-1)\right\}^2}\right]-
\E\left[e^{-\frac{t^2}{2}X_0^2}\right]\right|.
$$
The right-hand side of the previous equality being a $O(\frac1k)$ by 
the classical Berry-Esseen theorem (e.g. with the Wasserstein distance), the 
condition (B1) of Theorem \ref{thm2} is immediately satisfied.

We now turn to (B2).
Using the product formula (\ref{multiplication}) for multiple integrals, we can write
\[
\sum_{k=1}^{n}\frac{G_k^2-\E[G_k^2]}{k}=
I_{6}\left(\sum_{k=1}^{n} \frac{1}{k}g_k\widetilde{\otimes} g_k\right)
+9\,I_{4}\left(\sum_{k=1}^{n} \frac{1}{k}g_k\widetilde{\otimes}_1 g_k\right)
+18\,I_{2}\left(\sum_{k=1}^{n} \frac{1}{k}g_k\widetilde{\otimes}_2 g_k\right)
.
\]
We are thus left to check that, for $r=\{0,1,2\}$:
$$
\sum_{n\geq 2}\frac{1}{n\log^3n}\left\|
\sum_{k=1}^n \frac{1}{k}\,g_k\otimes_r g_k
\right\|_{\HH^{\otimes(2q-2r)}}^2<\infty.
$$
This is actually immediate, by using the following
straightforward computations (left to the reader):
\begin{eqnarray}
\langle g_k\otimes g_k,g_l\otimes g_l\rangle_{\HH^{\otimes 6}}&=&\frac{(k\wedge l)^2}{36\,kl},\quad
\langle g_k\otimes_1 g_k,g_l\otimes_1 g_l\rangle_{\HH^{\otimes 4}}=\frac{4\,k\wedge l+(k\wedge l)^2}{324\,kl}\notag\\
\langle g_k\otimes_2 g_k,g_l\otimes_2 g_l\rangle_{\HH^{\otimes 2}}&=&\frac{5\,k\wedge l}{324\,kl}.\label{est}
\end{eqnarray}

Summarizing, we have shown the following result.
\begin{prop}
Let $X_0,X_1,X_2\ldots\sim N(0,1)$ be a sequence of independent copies
defined on the same probability space, and consider
$
G_n=\frac{X_n}{\sqrt{2n}}\sum_{j=0}^{n-1} (X_j^2-1).
$
Then $(G_n)$ satisfies an ASLT with $\mu$ the distribution of $X_0X_1$.
\end{prop}

{\it Remark}. It is interesting to observe that the expressions (\ref{est})
are no longer valid if we replace $G_n$ by 
$
\widetilde{G}_n=\frac{X_0}{\sqrt{2n}}\sum_{j=1}^{n} (X_j^2-1).
$
Actually, the assumption (B2) of Theorem \ref{thm2} turns out to be not satisfied for $(\widetilde{G}_n)$. We believe that it is because $\widetilde{G_n}$ does not satisfy an ASLT, see also the discussion in the second paragraph following Definition \ref{defi1}.

\section{Proof of Corollary \ref{thm3}}\label{section 4}

Let $(G_n)$ be given by (\ref{gn}) with $H\in (1-\frac1{2q},1)$
and recall $V_n$ from (\ref{vn}).
It is a straightforward exercise to prove that ${\rm Var}(V_n)$ is equivalent to 
a constant times  $n^{2-2q(1-H)}$ as $n\to\infty$ (see, e.g. \cite{BNT} and references therein).
Thanks to Lemma \ref{reduction} and because the Hermite distribution $\mu_{H,q}$ admits a density (this latter fact is an immediate consequence of the main result of \cite{nourdinpoly} for instance), it is then equivalent to prove
an ASLT for $(G_n)$ or for ($\widehat{G}_n$) defined as
\[
\widehat{G}_n=n^{q(1-H)-1}\sum_{j=1}^{n}H_q(B_{j+1}-B_j).
\]
In the sequel, we are thus rather going to prove that $(\widehat{G}_n)$
satisfies an ASLT.

Our first observation is that $\widehat{G}_n$ takes the form of a multiple Wiener-It\^o integral with respect to $B$, after identifying this latter with an isonormal Gaussian process $X$ over the Hilbert space $\HH$ defined as the closure of the set of step functions with respect to the scalar product 
$\langle {\bf 1}_{[0,t]},{\bf 1}_{[0,s]}\rangle_\HH = \E[B_tB_s]$, $s,t\geq 0$, see \cite[Proposition 7.2.3]{n-p-book} for further details.
More precisely, we have
\[
\widehat{G}_n=I_q(g_n),\quad \mbox{where $g_n=n^{q(1-H)-1}\sum_{j=1}^{n}{\bf 1}_{[j,j+1]}^{\otimes q}$},
\]

This observation being made, we can now use Theorem \ref{thm2}, and check that conditions (B1) and (B2) therein are in order. In what follows, $\phi_k$ and $\phi_\infty$ denote the characteristic functions of $\widehat{G}_k$ and $\mu_{H,q}$ respectively. 

For (B1), it is a direct consequence of some estimates given in \cite{BN}.
More precisely,  using \cite[Proposition 3.1]{BN} one infers that, for any (fixed) $r>0$,
\begin{equation}\label{eq:B1}
\sup_{|t|\leq r}|\phi_k(t)-\phi_\infty(t)|\leq r d_W(\widehat{G}_k,\mu_{H,q})=O( k^{q(1-H)-\frac12}),
\end{equation}
where $d_W$ stands for the Wasserstein distance. The first inequality in (\ref{eq:B1}) is a direct consequence of the basic fact that $x \mapsto e^{itx}$ is a Lipschitz-continuous function with Lipschitz constant $\vert t \vert $. Also, recall for a real-valued random variables $F_1\sim \mu_1$ and $F_2\sim \mu_2$ that 
$$
d_W (F_1,F_2 )= d_W(\mu_1,\mu_2)
=\sup \left\{ \left \vert \int_\R hd\mu_1 - \int_\R hd\mu_2 \right \vert \, : 
h:\R \to \R \, \text{ s.t. } \, \Vert h \Vert_{\text{Lip} } \le 1  \right\}.
$$  Therefore, for some constant $C$, we obtain
\[
\sup_{|t|\leq r}\sum_{n\geq 2}\frac{1}{n\log^3n}\left|\sum_{k=1}^n\frac1k(\phi_k(t)-\phi_\infty(t))\right|^2
\leq 
C \sum_{n\geq 2}\frac{1}{n\log^3n}\left(\sum_{k=1}^n
k^{q(1-H)-\frac32}\right)^2
<\infty.
\]

Let us now turn to (B2).
Using the product formula (\ref{multiplication}) for multiple integrals, we can write
\[
\sum_{k=1}^{n}\frac{G_k^2}{k}=\sum_{r=0}^q r!\binom{q}{r}^2
I_{2q-2r}\left(\sum_{k=1}^{n} \frac{1}{k}g_k\widetilde{\otimes}_r g_k\right).
\]
Thus we are left to check that, for all $r\in\{0,\ldots,q-1\}$:
\begin{equation}\label{tocheck}
\sum_{n\geq 2}\frac{1}{n\log^3n}\left\|
\sum_{k=1}^n \frac{1}{k}\,g_k\otimes_r g_k
\right\|_{\HH^{\otimes(2q-2r)}}^2<\infty.
\end{equation}
Let us first focus on the case where $1\leq r\leq q-1$.
We can write, with $\rho(a)=\frac{1}{2}\big(|a+1|^{2H}+|a-1|^{2H}-2|a|^{2H}\big)$
 and $\rho_k(a)=|\rho(a)|{\bf 1}_{|a|\leq k}$:
\begin{eqnarray*}
&&\left\|
\sum_{k=1}^n \frac{1}{k}\,g_k\otimes_r g_k
\right\|_{\HH^{\otimes(2q-2r)}}^2=\sum_{k,l=1}^n \frac{1}{kl}\langle g_k\otimes_r g_k,g_l\otimes_r g_l\rangle_{\HH^{\otimes(2q-2r)}}\\
&\leq&2\sum_{1\leq l\leq k\leq n} (kl)^{2q(1-H)-3}\sum_{i,j=1}^{k}\sum_{i',j'=1}^{l}
\rho(i-j)^r\rho(i-i')^{q-r}
\rho(i'-j')^r\rho(j-j')^{q-r}\\
&\leq&2\sum_{1\leq l\leq k\leq n} (kl)^{2q(1-H)-3}\sum_{i=1}^{k}\sum_{j'=1}^{l}
\big(\rho_k^r * \rho_k^{q-r})(i-j')^2\\
&\leq&2\sum_{1\leq l\leq k\leq n} (kl)^{2q(1-H)-3}l
\|\rho_k^r * \rho_k^{q-r}\|^2_{\ell^2(\mathbb{Z})}.
\end{eqnarray*}
Using the Young inequality with $p=\frac{2q}{3r}$ and $p'=\frac{2q}{3(q-r)}$
(so that $\frac{1}p+\frac1{p'}=\frac32$), we obtain
\[
\|\rho_k^r * \rho_k^{q-r}\|^2_{\ell^2(\mathbb{Z})}
\leq \|\rho_k^r \|^2_{\ell^p(\mathbb{Z})}
\|\rho_k^{q-r} \|^2_{\ell^{p'}(\mathbb{Z})}
=  \left(
\sum_{|j|\leq k}|\rho(j)|^{\frac{2q}{3}}
\right)^3.
\]
But $\rho(j)\sim H(2H-1)|j|^{2H-2}$ and, because $H>1-\frac{1}{2q}>1-\frac{3}{4q}$, we deduce that
\[
\left(
\sum_{|j|\leq k}|\rho(j)|^{\frac{2q}{3}}
\right)^3=O( k^{3-4q(1-H)})\quad \mbox{as $k\to\infty$}.
\]
Thus, for $r\in\{1,\ldots,q-1\}$,
\begin{eqnarray}\notag
\left\|
\sum_{k=1}^n \frac{1}{k}\,g_k\otimes_r g_k
\right\|_{\HH^{\otimes(2q-2r)}}^2&=&O\left(\sum_{1\leq l\leq k\leq n} 
l^{2q(1-H)-2}k^{-2q(1-H)}\right)=O\left(\sum_{k=1}^n \frac1k\right) \\
&=& O(\log n). \label{csq1}
\end{eqnarray}
Actually, the previous estimate (\ref{csq1}) is also valid for $r=0$. Indeed, we have in this case
\begin{eqnarray}
&&\left\|
\sum_{k=1}^n \frac{1}{k}\,g_k\otimes_0 g_k
\right\|_{\HH^{\otimes(2q)}}^2=\sum_{k,l=1}^n (kl)^{2q(1-H)-3}\left(
\sum_{i=1}^k\sum_{j=1}^l \rho(i-j)^q
\right)^2\notag\\
&\leq& 2\sum_{1\leq l\leq k\leq n} (kl)^{2q(1-H)-3}\left(
\sum_{i=1}^k\sum_{j=1}^l |\rho(i-j)|^q
\right)^2\notag\\
&\leq& 2\sum_{1\leq l\leq k\leq n} (kl)^{2q(1-H)-3} l^2\left(
\sum_{|r|<k} |\rho(r)|^q
\right)^2\notag\\
&\leq&{\rm cst}\sum_{1\leq l\leq k\leq n} (kl)^{2q(1-H)-3}l^2 k^{2-4(1-H)q}
=O(\sum_{k=1}^n \frac{1}{k})=O(\log n),\label{csq2}
\end{eqnarray}
 where the last inequality used the fact that $\rho(r)\sim H(2H-1)r^{2H-2}$
as $|r|\to\infty$.
Finally,
 (\ref{csq1}) and (\ref{csq2}) together imply that 
(\ref{tocheck}) takes place for any $r\in\{0,\ldots,q-1\}$, which concludes the proof
of Theorem \ref{thm3}.
\qed

\section{Theorem \ref{thm2} implies Theorem \ref{thm1}}\label{section 5}

To conclude this paper, let us explain how to deduce Theorem \ref{thm1}
(taken from \cite{BNT})
from our Theorem \ref{thm2}, even if at first glance conditions (A1)-(A2) and (B1)-(B2) do not seem to be related to each others.

Fix $q\geq 2$, and let $(G_n)$ be a sequence of the form $G_n=I_q(g_n)$, 
$n\geq 1$. Suppose moreover that $\E[G_n^2]=q!\|g_n\|^2_{\HH^{\otimes q}}=1$ for all $n$.
In what follows, $c_q>0$ denotes a constant only depending on $q$, whose value can change from line to line.

Firstly, we deduce from \cite[Theorem 5.2.6]{n-p-book} (to get the second inequality) and \cite[(5.2.6)]{n-p-book} (to get the third inequality) that,
for all $t\in\R$,
\begin{eqnarray*}
\left| \E[e^{itG_k}]-e^{-t^2/2}\right|&\leq& |t|\,d_W(G_k,N(0,1))
\leq c_q\,|t|\,(\E[G_k^4]-3)\\
&\leq& c_q\,|t|\max_{1\leq r\leq q-1}\|g_k\widetilde{\otimes}_r g_k\|_{\HH^{\otimes(2q-2r)}}\\
&\leq&
c_q\,|t|\max_{1\leq r\leq q-1}\|g_k\otimes_r g_k\|_{\HH^{\otimes(2q-2r)}}.
\end{eqnarray*}

As a consequence,
\begin{eqnarray*}
&&\sup_{|t|\leq r}\sum_{n\geq 2}\frac{1}{n\log^3n}\left|\sum_{k=1}^n\frac1k(\E[e^{itG_k}]-e^{-t^2/2})\right|^2\\
&\leq& 2\,\sup_{|t|\leq r}
\sum_{n\geq 2}\frac{1}{n\log^2n}\left|\sum_{k=1}^n\frac1k(\E[e^{itG_k}]-e^{-t^2/2})\right|\\
&\leq& 2\,r\,c_q\,\max_{1\leq r\leq q-1}\sum_{n\geq 2}\frac{1}{n\log^2n}\sum_{k=1}^n\frac1k\,\|g_k\otimes_r g_k\|_{\HH^{\otimes(2q-2r)}},
\end{eqnarray*}
from which it comes that
(A1) implies (B1).

Secondly, one can write
\begin{eqnarray*}
\|g_k\otimes_r g_l\|_{\HH^{\otimes(2q-2r)}}^2
&=&\langle g_k\otimes_r g_k, g_l\otimes_r g_l\rangle_{\HH^{\otimes(2q-2r)}}\\
&\leq&\| g_k\otimes_r g_k\|_{\HH^{\otimes(2r)}} \|g_l\otimes_r g_l\|_{\HH^{\otimes(2q-2r)}}\\
&\leq&\| g_l\|^2_{\HH^{\otimes q}} \|g_k\otimes_r g_k\|_{\HH^{\otimes(2q-2r)}}
=q!\|g_k\otimes_r g_k\|_{\HH^{\otimes(2q-2r)}}.
\end{eqnarray*}
Moreover, it has been shown in \cite[Proof of Theorem 4.3]{NR}
that 
\begin{eqnarray*}
&&{\bf Cov}(G_k^2,G_l^2)-2\left(\E[G_kG_l]\right)^2\\
&=&q!^2\sum_{r=1}^{q-1}\binom{q}{r}^2\|g_k\otimes_r g_l\|_{\HH^{\otimes(2q-2r)}}^2+\sum_{r=1}^{q-1}r!^2\binom{q}{r}^4(2q-2r)!
\|g_k\widetilde{\otimes}_r g_l\|_{\HH^{\otimes(2q-2r)}}^2.
\end{eqnarray*}
We deduce that
\begin{eqnarray*}
0\leq {\bf Cov}(G_k^2,G_l^2)&\leq &2\left(\E[G_kG_l]\right)^2 +c_q 
\max_{1\leq r\leq q-1} \|g_k\otimes_r g_k\|_{\HH^{\otimes(2q-2r)}}\\
&\leq &2\left|\E[G_kG_l]\right| +c_q 
\max_{1\leq r\leq q-1} \|g_k\otimes_r g_k\|_{\HH^{\otimes(2q-2r)}},
\end{eqnarray*}
implying in turn that
\begin{eqnarray*}
&&
\sum_{n\geq 2}\frac{1}{n\log^3n}{\bf Var}\left(\sum_{k=1}^n \frac{G_k^2}{k}\right)=\sum_{n\geq 2}\frac{1}{n\log^3n}\sum_{k,l=1}^n \frac1{kl}{\bf Cov}\left(G_k^2,G_l^2\right)\\
&\leq&2\sum_{n\geq 2}\frac{1}{n\log^3n}\sum_{k,l=1}^n \frac{\left|\E[G_kG_l]\right|}{kl}+c_q\,\max_{1\leq r\leq q-1}
\sum_{n\geq 2}\frac{1}{n\log^2n}\sum_{k=1}^n \frac1{k}\|g_k\otimes_r g_k\|_{\HH^{\otimes(2q-2r)}},
\end{eqnarray*}
from which it follows that (A1) and (A2) imply (B2).\qed

\bigskip

{\bf Acknowledgment}. 
We thank two anonymous referees for their useful comments that led to an improvement in the presentation of this paper.

\end{document}